\documentclass[12pt]{amsart}
\usepackage{amsfonts}
\usepackage{amssymb}
\usepackage{amsmath,mathrsfs}
\newtheorem{theo}{Theorem}
\newtheorem{pro}{Proposition}
\newtheorem{defi}{Definition}
\newtheorem{lem}{Lemma}

\newtheorem*{rem}{Remark}
\usepackage{a4}
\usepackage{url}
\usepackage{color}

\newcommand{\bmultg}{\!\begin{array}{c} {\scriptstyle\times}
    \\[-12pt]\cup\end{array}\!}

\DeclareMathAlphabet{\eusm}{U}{}{}{}  
\SetMathAlphabet\eusm{normal}{U}{eus}{m}{n}
\SetMathAlphabet\eusm{bold}{U}{eus}{b}{n}

\DeclareMathAlphabet{\eufrak}{U}{}{}{}  
\SetMathAlphabet\eufrak{normal}{U}{euf}{m}{n}
\SetMathAlphabet\eufrak{bold}{U}{euf}{b}{n}

\theoremstyle{definition}

\theoremstyle{remark}

\numberwithin{equation}{section}

\makeatletter
\@namedef{subjclassname@2010}{%
  \textup{2010} Mathematics Subject Classification}
\makeatother

\begin{document}

\title{Monotone and Boolean unitary Brownian motions}

\author[T. Hamdi]{Tarek Hamdi}
\address{Universit\'e de Tunis El-Manar, Laboratoire d'Analyse Mathématiques et applications, LR 11ES11.}
\email{tarek.hamdi@ipest.rnu.tn}

\begin{abstract}
The additive monotone (resp. boolean) unitary Brownian motion is a non-commutative stochastic process with monotone (resp.
boolean) independent and stationary increments which are distributed according to the arcsine law (resp. Bernoulli law) .
We introduce the monotone and boolean unitary Brownian motions and we derive a closed formula for their associated  moments. This provides a description of their spectral measures. We prove that, in the monotone case, the multiplicative analog of the arcsine distribution is absolutely continuous with respect to the Haar measure on the unit circle, whereas in the boolean case the multiplicative analog of the Bernoulli distribution is discrete. Finally, we use quantum stochastic calculus
to provide a realization of these processes as the stochastic exponential of the correspending additive Brownian motions.
\end{abstract}

\keywords{Convolution semi-groups, monotone and
boolean independence, unitary Brownian motion, boson Fock space, boolean Fock space, quantum stochastic calculus.}

\subjclass[2010]{Primary 60J65, 46L51, 46L53; Secondary 65C30}
\maketitle

\section*{Introduction}
\label{sec-intro}

In non-commutative probability theory, there exist several natural notions of independence. 
The most famous ones are the usual independence and the free independence.
Other interesting examples are monotone and boolean independence. 
These allow to define new convolutions for probability measures.
The monotone convolutions on the unit circle and the positive half-line were
introduced by Bercovici in \cite{bercovici04}, see also \cite{franz05a}.
But the additive monotone Brownian motion and the monotone Fock space where first
studied by Muraki, see \cite{muraki97}. The monotone stochastic calculus can also be realized on
the symmetric Fock space, see \cite{franz03b}.
Bercovici studied also the boolean convolution for probability measures on the positive half-line, but it is not always defined. 
 For probability
measures on the unit circle, the boolean convolution is well defined. It was introduced by Franz in \cite{franz04}, see also \cite{bercovici04b}.
 The boolean stochastic calculus has been studied by
Ben Ghorbal and Sch\"urmann, see \cite{benghorbel+surmann04}.

The aim of this paper is to point out several connections between classical, free Brownian motions and their counterparts in the monotone and boolean cases.
We shall consider two kinds of unitary Brownian motions, the monotone and boolean one. Both are non-commutative unitary processes, that is families of 
non-commutative random variables which are unitary and are characterized as having independent increments, distributed according
 to a multiplicative convolution semi-group of measures depending on the notion of independence that we use (monotone or boolean).
It has been shown that the additive monotone Brownian motion is distributed according to the arcsine law (cf.\ \cite{muraki97})
\[
 \frac{1}{\pi\sqrt{2t-x^2}} \mathbf{1}_{(-\sqrt{2t},\sqrt{2t})} {\rm
  d}x,\quad t>0.
\]
Whereas the additive Boolean Brownian motion is Bernoulli distributed (cf.\ \cite{speicher+woroudi93})
\begin{equation*}
\frac{1}{2}(\delta_{-\sqrt{t}} + \delta_{\sqrt{t}}),\quad t>0.
\end{equation*}
So the distributions of the monotone and boolean unitary Brownian motions are a kind of
multiplicative analog of these two distributions. 
We shall next consider the Fock space realization of both unitary Brownian motion by solving the exponential stochastic differential equation
$$dU_t=i dX_t\, U_t-\frac{1}{2}U_t\, dt,\quad U_0=I$$
where $(X_t)_t$ is the Fock space realization of the corresponding additive Brownian motion. Note that this equation is formally the same as the one used to 
construct the free unitary Brownian motion (cf. \cite{biane97}).

This paper is organized as follows. In part one, we first introduce the monotone unitary Brownian motion.
We then derive a closed formula for its moments through Legendre polynomials and supply a full description of its distribution, say $\mu_t$.
In particular, $\mu_t$ is compactly-supported and  absolutely continuous with respect to the Haar
measure on the unit circle.
Moreover, it's support is inside an interval of the unit circle, centered at 1, and spreads as $t$ increases. It only becomes the whole circle asymptotically as $t$ goes to infinity.
We close the first part by giving the boson Fock space realization of the monotone unitary Brownian motion. 
The second part is concerned with boolean unitary Brownian motion. We first introduce this process and we
derive an explicit formula for its moments which involves Laguerre polynomials. 
It turns out that the  distribution of the boolean unitary Brownian motion, say $\nu_t$, is discrete. Its support is a countably set of points  on the unit circle which accumulate at 1.
Finally, we recall some basic facts from \cite{benghorbel+surmann04} and we give a Fock space realization of the boolean unitary Brownian motion.

\section{\bf Monotone case}
\subsection{Multiplicative monotone convolution}

Let $\mu_1,\mu_2$ be probability measures on the unit circle and set
\[
\psi_{\mu_i} = \int_{\mathbb{T}} \frac{zx}{1-zx}{\rm d}\mu_i(x), \qquad
K_{\mu_i}(z) = \frac{\psi_{\mu_i}(z)}{1+\psi_{\mu_i}(z)},
\]
then the multiplicative monotone convolution of $\nu=\mu_1 \triangleright \mu_2$ is uniquely
determined  by
\[
K_\nu(z) = K_{\mu_1}\big(K_{\mu_2}(z)\big).
\]
The monotone convolution semigroups of probability measures on the unit circle are
given by solutions of the differential equation
\[
\frac{{\rm d}K_t}{{\rm d}t} (z) = - K_t(z)u\big(K_t(z)\big)
\]
with initial condition $K_0(z)=z$, where $u$ is a holomorphic function on the
unit disk with non-negative real part. Therefore $u$ has a Herglotz
representation
\[
u(z) = ib + \int \frac{x+z}{x-z}{\rm d}\rho(x)
\]
with $b$ a real number and $\rho$ a finite measure supported on the unit
circle.
\subsection{Monotone unitary Brownian motion}
The distribution $\mu_t$  of the unitary monotone Brownian motion   
correspond to the case $b=0$ and $\rho=\frac{1}{2}\delta_1$, i.e. the transform
$(K_{\mu_t})_{t\ge 0}$ is given by the solution of the differential equation
\begin{equation}\label{ode}
2\frac{{\rm d}K_{\mu_t}}{{\rm d}t} (z) = - K_{\mu_t}(z)\frac{1+K_{\mu_t}(z)}{1-K_{\mu_t}(z)},
\end{equation}
with initial condition $K_{\mu_0}(z)=z$.
If we write $Z_t=K_{\mu_t}$, then the differential equation \eqref{ode} becomes
\[
2\dot{Z}_t = \frac{Z_t^2 + Z_t}{Z_t-1}
\]
or
\[
\left(\frac{2}{Z_t+1}- \frac{1}{Z_t}\right) = \frac{1}{2},
\]
which can be integrated to
\[
\frac{t}{2}=\int_z^{Z_t} \left(\frac{2}{u+1}- \frac{1}{u}\right){\rm d}u = 2
\ln\frac{Z_t+1}{z+1} - \ln \frac{Z_t}{z}.
\]
Exponentiating yields
\[
\frac{(Z_t+1)^2}{Z_t} = e^{t/2} \frac{(z+1)^2}{z}.
\]
We set $w=e^{t/2} \frac{(z+1)^2}{z}$ and solve the quadratic equation
\begin{equation}\label{equ}
Z_t^2 +(2-w)Z_t+1=0
\end{equation}
for $Z_t$. This gives
\begin{equation*}
Z_t = \frac{1}{2}\left(w-2\pm\sqrt{w^2-4w}\right)
\end{equation*}
where the sign $\pm$ is such that $Z_t$ is inside the unit disk and
$\lim_{t\searrow 0}=z$. 
But since the only bounded solution of (\ref{equ}) is what corresponds to a minus sign, then we get
\begin{equation}\label{K}
Z_t= \frac{1}{2}\left( w-2-\sqrt{w(w-4)}\right).
\end{equation} 
Now the moment generating function
\[
\psi_{\mu_t}(z) := \sum_{k=1}^\infty z^k\int_{\mathbb{T}} x^k {\rm d}\mu_t(x) = \int_{\mathbb{T}} \frac{xz}{1-xz}{\rm d}\mu_t(x)
\]
of the distribution of unitary monotone Brownian motion can be recovered
as
\[
\psi_{\mu_t} = \frac{Z_t}{1-Z_t}.
\]
Note in passing that the family $(K_t)_{t\ge 0}$ defined by $K_t(z)=Z_t|_{Z_0=z}$
is a family of self maps of the unit disk which fixes the origin, i.e., $K_t(0)=0$.
Besides it forms a continuous composition semi-group and it is a
general feature that these composition semi-group have such a form.
\begin{rem}
Note that the equality (\ref{K}) can also be written as
\[
Z_t(z) = \varphi^{-1}\left(e^{t/2} \varphi(z)\right)
\]
where $\varphi$ is the conformal bijection from $\mathbb{D}$ to $\mathbb{C}\setminus[0,4]$ given by $\varphi(z)=(z+1)^2/z$.
\end{rem}

\subsection{Moments of $\mu_t$}
\begin{pro}\label{Mom1}
 For every $t>0$, one has
\begin{equation*}
 \int_{\mathbb{T}} x^n {\rm d}\mu_t(x)=\frac{1}{2}\left(P_n(2e^{-t/2}-1)+P_{n-1}(2e^{-t/2}-1)\right),\quad n \geq 1
\end{equation*}
 where $P_n$ is the $n$-th Legendre polynomial.
\end{pro}
{\bf Proof}
Since we have
\begin{align*}
 Z_t(z)= \frac{1}{2}\left( w-2-\sqrt{w(w-4)}\right)\quad z\in \mathbb{D}
\end{align*}
 with $w=e^{t/2}\varphi(z)$.
Then we obtain
\begin{align*}
 \psi_{\mu_t}(z)=\frac{Z_t}{1-Z_t}&=-\frac{w-\sqrt{w(w-4)}-2}{w-\sqrt{w(w-4)}-4}
\\&=-1+\frac{-2}{w-\sqrt{w(w-4)}-4}.
\end{align*}
Writing
\begin{align*}
 \frac{1}{w-\sqrt{w(w-4)}-4}=\frac{-1}{4}\left( 1+\sqrt{\frac{w}{w-4}}\right) ,
\end{align*}
 we get
\begin{align*}
\psi_{\mu_t}(z)&=-\frac{1}{2}+\frac{z+1}{2\sqrt{1-2z(2e^{-t/2}-1)+z^2}}.
\end{align*}
Now using the generating function for Legendre polynomial 
\begin{align*}
\frac{1}{\sqrt{1-2xz+z^2}} = \sum_{n\geq 0} P_n(x) z^n,
\end{align*}
we obtain the following expansion
\begin{align*}
 \psi_{\mu_t}(z)&=-\frac{1}{2}+\frac{z+1}{2}\sum_{n\geq0}P_n(2e^{-t/2}-1) z^n
\\&= -\frac{1}{2}+\frac{P_0(2e^{-t/2}-1)}{2}+\frac{1}{2}\sum_{n\geq1}(P_n(2e^{-t/2}-1)+P_{n-1}(2e^{-t/2}-1)) z^n
\\&= \frac{1}{2}\sum_{n\geq1}(P_n(2e^{-t/2}-1)+P_{n-1}(2e^{-t/2}-1)) z^n.
\end{align*}
 $\hfill\square$
\begin{rem}
The moments of the monotone unitary BM are real. This is due to the fact that the distribution
of the unitary monotone BM is symmetric, i.e., $U_t$ and
$U^*_t$ have the same distribution.
\end{rem}

\subsection{Description of $\mu_t$}

\begin{theo}
 For any $t>0$, the measure $\mu_t$ is absolutely continuous with respect to the Haar measure on $\mathbb{T}$. Its support is equal to the interval
\begin{equation*}
 I_t\triangleq \left\lbrace e^{i\theta}\left| -\arccos (2e^{-t/2}-1)\leq \theta \leq \arccos (2e^{-t/2}-1)\right. \right\rbrace 
\end{equation*}
where $\arccos\in[0,\pi]$. The density is positive on the interior of $I_t$ and is equal to
\begin{equation*}
 \frac{\sqrt{2}\cos (\theta/2)}{\sqrt{\cos (\theta)-(2e^{-t/2}-1)}}1_{\{\cos (\theta)>2e^{-t/2}-1\}},\quad \theta\in (-\pi,\pi).
\end{equation*}
\end{theo}

{ \bf Proof}
Let $H_{\mu_t}$ be the Herglotz transform of $\mu_t$. We have
\begin{align*}
H_{\mu_t}(z)&=1+2\psi_{\mu_t}(z)=\frac{z+1}{\sqrt{1-2z(2e^{-t/2}-1)+z^2}},\quad |z|<1.
\end{align*}
Since $z\mapsto 1-2z(2e^{-t/2}-1)+z^2$ does not take negative value and its roots $a_t, \overline{a_t}$ lie on $\mathbb{T}$ with
$
a_t=e^{i\arccos (2e^{-t/2}-1)}.
$
Then $H_{\mu_t}$ is analytic in the open disc and it extends 
continuously to $\mathbb{T} \setminus \{a_t,\overline{a_t}\}$. 
Furthermore for $\theta_t=\arccos (2e^{-t/2}-1)$, we have 
\begin{align*}
\lim_{r\rightarrow 1,z=re^{\pm i \theta_t}}(z-e^{\pm i \theta_t})H_{\mu_t}(z)= 0.
\end{align*}
 In fact, for any $z=re^{\pm i \theta_t}$, we write
\begin{align*}
(z-e^{\pm i \theta_t})H_{\mu_t}(z)&= (z-e^{\pm i \theta_t}) \frac{z+1}{\sqrt{(z-a_1)(z-a_2)}}
\\&= (z-e^{\pm i \theta_t}) \frac{z+1}{\sqrt{(z-e^{ i \theta_t})(z-e^{- i \theta_t})}}
\\&=\frac{(r-1)e^{\pm i \theta_t}(e^{\pm i \theta_t}+1)}{\sqrt{r-1}\sqrt{re^{\pm 2i \theta_t}-1}}
\\&=\frac{\sqrt{r-1}e^{\pm i \theta_t}(e^{\pm i \theta_t}+1)}{\sqrt{re^{\pm 2i \theta_t}-1}}.
\end{align*}
Then $\mu_t$ is absolutely continuous with respect to the Haar measure on the unit circle, with density given by the Poisson kernel
$$P_{\mu_t}(e^{i\theta})=\Re\left[H_{\mu_t}(e^{i\theta}) \right],\quad \theta\in (-\pi,\pi).$$
Now we need to compute $\Re\left[H_{\mu_t}(w) \right]$ for $w\in \mathbb{T}$.  We have
\begin{align*}
H_{\mu_t}(w)&=\frac{w+1}{\sqrt{1-2w(2e^{-t/2}-1)+w^2}}
\\&=\frac{w+1}{\sqrt{w}\sqrt{w+\overline{w}-2(2e^{-t/2}-1)}}1_{\{(w+\overline{w})-2(2e^{-t/2}-1)>0\}}
\\&=\frac{(w+1)\left(\sqrt{(2+w+\overline{w})/4}-i\sqrt{(2-w-\overline{w})/4}\right)}{\sqrt{w+\overline{w}-2(2e^{-t/2}-1)}}1_{\{(w+\overline{w})-2(2e^{-t/2}-1)>0\}}
\\&=\frac{(2+w+\overline{w}+(w-\overline{w}))\left(\sqrt{(2+w+\overline{w})/4}-i\sqrt{(2-w-\overline{w})/4}\right)}{2\sqrt{w+\overline{w}-2(2e^{-t/2}-1)}}1_{\{(w+\overline{w})-2(2e^{-t/2}-1)>0\}}.
\end{align*}
Hence we get
\begin{align*}
 P_{\mu_t}(e^{i\theta})= &=\frac{\sqrt{2}\cos (\theta/2)}{\sqrt{\cos (\theta)-(2e^{-t/2}-1)}} 1_{\{\cos (\theta)>2e^{-t/2}-1\}}.
\end{align*}

$\hfill\square$
\begin{rem}
Note that for classical Brownian motion the support is immediately the whole circle and
in the free case, it is an interval smaller then the circle for $t<4$ and the whole circle for $t\ge 4$. Whereas,
for the monotone case it remains smaller than the circle for all $t$, it only becomes the whole circle asymptotically as $t$ goes to infinity.
\end{rem}
\subsection{Quantum stochastic calculus and monotone unitary Brownian motion}
\subsubsection{\bf Quantum stochastic calculus}
We shall recall some basic facts from quantum stochastic calculus in a boson Fock space as expounded in \cite{parthasarathy99}.
Let $H:=\Gamma(L^2(\mathbb{R}_+))$ be the boson Fock space with one degree of freedom and let $A^*,A$ be respectively the creation and annihilation operators
on $H$. For each $t\in\mathbb{R}_+$ we have the following identification
$$\Gamma(L^2(\mathbb{R}_+))\cong \Gamma(L^2([0,t]))\otimes \Gamma(L^2([t,\infty))).$$
Let $\Omega_{t]}$ be the Fock vacuum vector in $H_{t]}:=\Gamma(L^2([0,t]))$ and let $P_0(t)$ be the projection on the subspace 
$\Omega_{t]}\otimes H_{[t}$ where $H_{[t}:=\Gamma(L^2([t,\infty)))$. According to \cite[Theorem 2.1]{parthasarathy99}, the quantum stochastic integrals
$$L_t=\int_0^tP_0(s)dA_s,\quad L_t^*=\int_0^tP_0(s)dA_s^*$$
are defined on the whole space $H$ as bounded operators, adjoint to each other and satisfy
$$L_t\,L_t^*=\int_0^tP_0(s)ds.$$
Now, let us denote $B_t=A_t+A_t^*$ and $X_t=L_t+L_t^*$. We recall that the family of operators $(X_t)_{t\in\mathbb{R}_+}$ on $\Gamma(H)$ is the non-commutative
arcsine Brownian motion for the vacuum expectation 
$\phi(.)=\langle e(0),e(0)\rangle$ on $B(\Gamma(H))$, (See \cite{parthasarathy99}).

\subsubsection{\bf Boson Fock space realization of monotone unitary Brownian motion}
We shall now use quantum stochastic calculus on a boson Fock space in order to give a construction of monotone unitary Brownian motion as the 
``stochastic exponential'' of the monotone additive one. For this we shall solve a stochastic differential equation.
According to \cite{hudson+parthasaraty84} there exists a unique family of unitary operators $(U_t)_{t\in\mathbb{R}_+}$, solution to the stochastic differential equation
\begin{equation}\label{QSDE}
 dU_t=(i\,dX_t-\frac{1}{2}P_0(t)dt)U_t,\quad U_0=1.
\end{equation}
>From this equation we deduce,
\begin{lem}\label{EQU3}
 For all $n\geq1$ and $t\geq 0$ one has
\begin{align*}
 dU_t^n=i\sum_{l=1}^nU_t^{n-l}P_0U_t^ldB_t-\frac{1}{2}\sum_{l=1}^nU_t^{n-l}P_0U_t^ldt-\sum_{k=2}^n\sum_{l=1}^{k-1}U_t^{n-k}P_0U_t^{k-l}P_0U_t^ldt.
\end{align*}
\end{lem}
{ \bf Proof}
We prove the result by induction on $n$. For $n=1$ this reduce to the equation for $U$. Assume the result holds for $n$, then by quantum Ito's formula
\begin{align*}
 dU_t^{n+1}=&U_t\,dU_t^n+dU_t\,U_t^n+dU_t\,dU_t^n
\\=&i\sum_{l=1}^nU_t^{n-l+1}P_0U_t^{l}dB_t-\frac{1}{2}\sum_{l=1}^nU_t^{n-l+1}P_0U_t^{l}dt-\sum_{k=2}^n\sum_{l=1}^{k-1}U_t^{n-k+1}P_0U_t^{k-l}P_0U_t^{l}dt
\\&+(i\,P_0dB_t-\frac{1}{2}P_0dt)U_t^{n+1}-\sum_{l=1}^nP_0U_t^{n-k+1}P_0U_t^{l}dt
\\=&i\sum_{l=1}^{n+1}U_t^{n+1-l}P_0U_t^{l}dB_t-\frac{1}{2}\sum_{l=1}^{n+1}U_t^{n+1-l}P_0U_t^{l}dt-\sum_{k=2}^{n+1}\sum_{l=1}^{k-1}U_t^{n+1-k}P_0U_t^{k-l}P_0U_t^{l}dt.
\end{align*}
$\hfill\square$

We shall now prove that
\begin{theo}\label{theo1}
 The process $(U_t)_{t\in\mathbb{R}_+}$ is a monotone unitary Brownian motion.
\end{theo}
{ \bf Proof} 
Form (\ref{QSDE}), we have
\begin{equation*}
 U_t-U_s=i\int_s^tdX_rU_r-\frac{1}{2}\int_s^tP_{0}(r)U_rdr,
\end{equation*}
so that,
\begin{equation}\label{increm}
 U_tU_s^{\star}=1+i\int_s^tP_{0}(r)U_rU_s^{\star}dB_r-\frac{1}{2}\int_s^tP_{0}(r)U_rU_s^{\star}dr.
\end{equation}
But the kernels of the projections $P_0(s)$ increase with $s$, this says that $P_0(r)=P_0(r)P_0(s)$ for any $s\leq r$. Consequently from (\ref{increm}), we get
$$U_tU_s^{\star}-1=P_0(s)(U_tU_s^{\star}-1)$$
and so for all $0\leq s_1<t_1<...<t_{n-1}< s_n<t_n<\infty$ and $k_1,...,k_n\in\mathbb{N}$, we get 
\begin{align*}
 \phi\left((U_{t_1}U_{s_1}^{\star}-1)^{k_1}...(U_{t_n}U_{s_n}^{\star}-1)^{k_n}\right)&=\phi\left(P_0(s_1)(U_{t_1}U_{s_1}^{\star}-1)^{k_1}...P_0(s_n)(U_{t_n}U_{s_n}^{\star}-1)^{k_n}\right)
\\&=\prod_{i=1}^n\phi\left( (U_{t_i}U_{s_i}^{\star}-1)^{k_i}\right) .
\end{align*}
Which proves that the increments of the process $(U_t)_{t\in\mathbb{R}_+}$ are monotonically independent with respect to the state $\phi$.
Now the equality (\ref{increm}) implies that the process $t\mapsto U_tU_s^{\star}$ is obtained by solving the same equation as $U$ but with respect to the monotone Brownian motion $(X_t-X_s)_{t\geq s}$. 
This says that the increments of $(U_t)_{t\in\mathbb{R}_+}$  have a stationary distribution. So it is enough 
 to prove that the distribution of $U_t$ is $\mu_t$ for all $t\geq 0$. To this end, we take the state $\phi$ in both sides of the equation in Lemma \ref{EQU3}. 
Then we get
\begin{align*}
 \frac{d}{dt}\phi(U_t^n)=-\frac{1}{2}\sum_{l=1}^n\phi(U_t^{n-l})\phi(U_t^{l})-\sum_{k=2}^n\sum_{l=1}^{k-1}\phi(U_t^{n-k})\phi(U_t^{k-l})\phi(U_t^{l}).
\end{align*}
 Now, introducing the generating function
\begin{equation*}
 \rho(t,z)=\sum_{n=1}^{\infty}\phi(U_t^{n})z^n,\quad |z|<1
\end{equation*}
we get the differential equation
\begin{align*}
\partial_t \rho(t,z)&=-\frac{1}{2}\rho(t,z)(1+\rho(t,z))-\rho(t,z)^2(1+\rho(t,z)),
\end{align*}
or equivalently
\begin{equation*}
 \partial_t \rho(t,z)=-\frac{1}{2}\rho(t,z)(1+\rho(t,z))(1+2\rho(t,z)).
\end{equation*}
After integrating and taking into account $\rho(0,z)=z/(1-z)$, we get
\begin{align*}
\frac{\rho(t,z)(1+\rho(t,z))}{(1+2\rho(t,z))^2}  =\frac{ze^{-\frac{t}{2}}}{(1+z)^2}=\frac{1}{\varphi(Z_t(z))}.
\end{align*}
Hence $\rho(t,.)=\psi_{\mu_t}$ for any $t\geq0$ and the distribution of $U_t$ is $\mu_t$. This ends the proof of Theorem \ref{theo1}.
$\hfill\square$

\section{\bf Boolean case}
\subsection{Multiplicative boolean convolution}

Let $\mu_1,\mu_2$ be probability measures on the unit circle and set
\[
\psi_{\mu_i} = \int_{\mathbb{T}} \frac{zx}{1-zx}{\rm d}\mu_i(x), \qquad
F_{\mu_i}(z) = \frac{1}{z} \frac{\psi_{\mu_i}(z)}{1+\psi_{\mu_i}(z)},
\]
then the multiplicative boolean convolution $\nu=\mu_1 \bmultg \mu_2$ is
uniquely determined  by
\[
F_\nu(z) = F_{\mu_1}(z)F_{\mu_2}(z).
\]

For infinitely divisible measures the transform is of the form
\[
F(z) = \exp\big(u(z)\big),
\]
where $u$ is a holomorphic function on the unit disk with non-positive real
part and therefore has a Herglotz representation of the form
\[
u(z) = ib - \int \frac{x+z}{x-z}{\rm d}\rho(x)
\]
with $b$ a real number and $\rho$ a finite measure supported on the unit
circle, cf.\ \cite[Remark 3.7]{franz04}.
\subsection{Boolean unitary Brownian motion}

The distribution $\nu_t$ of the boolean Brownian motion on the unit
circle corresponds again to the case $b=0$, $\rho=\frac{1}{2}\delta_1$, i.e.\ to the transforms
\[
F_t(z) = \exp\left(  \frac{t(z+1)}{2(z-1)}\right),\quad t\ge 0.
\]
In this case, one gets 
\begin{equation*}
\psi_{\nu_t}(z) = \frac{z\exp\left(  \frac{t(z+1)}{2(z-1)}\right)}{1-z\exp\left(  \frac{t(z+1)}{2(z-1)}\right)}
\end{equation*}
and in a small neighborhood of the origin 
\begin{equation*}
\psi_{\nu_t}(z) = \sum_{k \geq 1}z^k\exp\left(k \frac{t(z+1)}{2(z-1)}\right) = \sum_{k \geq 1}(ze^{-\frac{t}{2}})^k\exp\left(kt \frac{z}{z-1}\right) .
\end{equation*}

\subsection{Moments of $\nu_t$}
\begin{pro}\label{Mom2} 
 For any $t>0$ and any $n \geq 1$, one has
\begin{equation*}
 \int_{\mathbb{T}} x^n {\rm d}\nu_t(x)=\sum_{k=1}^nL_{n-k}^{(1)}(kt)e^{-kt/2} - 2\sum_{k=1}^{n-1}L_{n-k-1}^{(1)}(kt)e^{-kt/2} + \sum_{k=1}^{n-2}L_{n-k-2}^{(1)}(kt)e^{-kt/2}
\end{equation*}
where the sums in the RHS are taken to be empty when $k<1$ and $L^{(1)}_n$ is the $n$-th Laguerre polynomial of index 1.
\end{pro}
{\bf Proof}
We use the generating function for Laguerre polynomials 
\begin{equation*}
\exp\left(kt \frac{z}{z-1}\right) = (1-z)^{\alpha+1}\sum_{n \geq 0}L_n^{(\alpha)}(kt)z^n 
\end{equation*}
where $\alpha > -1$. We can choose $\alpha=1$ and get after inverting the order of summation 
\begin{align*}
\psi_{\nu_t}(z) &= (1-z)^2\sum_{k \geq 1}e^{-kt/2}\sum_{n \geq k}L_{n-k}^{(1)}(kt)z^{n} 
\\& =  (1-z)^2\sum_{n \geq 1} \left(\sum_{k=1}^nL_{n-k}^{(1)}(kt)e^{-kt/2}\right)z^n.
\end{align*}
 Identifying the coefficients of $z^n$, we get the result.
  $\hfill\square$

\begin{rem}
The moments are again real, because the law of the unitary
  boolean BM is also symmetric.
 
\end{rem}

\subsection{Description of $\nu_t$}

 Observe that the Herglotz transform $k_t$ of $\nu_t$ is given by
\begin{equation*}
 k_t(z)=\int_{\mathbb{T}}
\frac{\xi+z}{\xi-z}d\nu_t(\xi)=1+2\psi_{\nu_t}(z)=\frac{1+\theta_t(z)}{1-\theta_t(z)},\quad t>0, z\in \mathbb{D}
\end{equation*}
with
$$\theta_t(z)=z\exp\left(  \frac{t(z+1)}{2(z-1)}\right).$$
 For $t> 0$, let $\Sigma_t$ denote the subset of $\mathbb{T}$ consisting of the solutions to the equation $\theta_t(\zeta)=1$.
We shall need the following lemma.
\begin{lem}\label{discret}
$\Sigma_t$ is a discrete set of points which accumulate at 1.
\end{lem}
{ \bf Proof}
Fix $t>0$ and let $\zeta=x+iy \in\Sigma_t$. This says that
\begin{align*}
\exp\left( -\frac{it y}{2(1-x)}\right) = x-iy.
\end{align*}
By a simple identification of real and imaginary parts we obtain
\begin{align*}
  \begin{cases}
                    \frac{ t y}{2(1-x)}  = \epsilon \arccos (x)+2k\pi,\ \epsilon\in\{\pm 1\}\ \text{and}\ k\in\mathbb{Z}\\
\sin \left(\frac{ t y}{2(1-x)} \right) = y
                   \end{cases}
\end{align*}
which in turn leads to
\begin{align*}
 \begin{cases}
                    \frac{t}{2}\sqrt{\frac{ 1+x}{1-x} } - \arccos (x) = 2k\epsilon\pi,\ \epsilon\in\{\pm 1\}\ \text{and}\ k\in\mathbb{Z}\\
  y = \epsilon \sqrt{1-x^2}
                   \end{cases}.
\end{align*}
Consider now the function $g_t(x)= \frac{t}{2}\sqrt{\frac{ 1+x}{1-x} } - \arccos (x)$ with $x\in[-1,1)$. Then we have
\begin{equation*}
 \partial_x g_t(x)=\frac{t+2-2x}{2(1-x)\sqrt{1-x^2}}\geq0
\end{equation*}
and so \\

\centerline{$
\begin{array}{|c|ccccr|}
\hline
x     & -1   & & x_0(t)  &    &1  \\
\hline
\partial_x g_t(x) &  &    & +  &  & \\
\hline &&  &&   & \\       
g_t(x) &-\pi  &\nearrow  & 0 &\nearrow  & +\infty \\   
&  & & & &  \\         
\hline
\end{array}
$
}
\vspace{0.7cm}

Consequently for any $k\in\mathbb{Z}$, the equation $g_t(x)=2k\epsilon \pi$ has a unique solution $x_k(t)$ if $ \epsilon k\geq0$ and no solution otherwise.
Hence
\begin{equation*}
\Sigma_t=\left\lbrace x_n(t)+ i\sqrt{1-x_n(t)^2},\ n\in\mathbb{N}\right\rbrace\cup \left\lbrace x_{-n}(t)- i\sqrt{1-x_{-n}(t)^2},\ n\in\mathbb{N}\right\rbrace.
\end{equation*}
Since $g_t$ is injective and we have
\begin{equation*}
 g_t(x_{-n}(t))=2(-n)(-1)\pi=2n\pi=g_t(x_n(t)),
\end{equation*}
we obtain $x_n(t)=x_{-n}(t)$ then

\begin{equation*}
\Sigma_t=\left\lbrace \zeta_{\pm n}(t)=x_n(t)\pm i\sqrt{1-x_n(t)^2}:\ g_t(x_n(t))=2n\pi,n\in\mathbb{N}\right\rbrace.
\end{equation*}
It remains to show that $(\zeta_n(t))_{n\in\mathbb{Z}}$ accumulate at 1.
Equivalently, we shall show that
$\lim_{n \rightarrow +\infty} x_n(t)=1$. 
Notice that the sequence $(x_n(t))_{n\in\mathbb{N}}$ is bounded. Moreover, for $n\leq m$, we have
$$
g_t(x_{n}(t))=2\pi n\leq 2\pi m=g_t(x_{m}(t)).
$$
Thus $(x_n(t))_{n\in\mathbb{N}}$ is increasing and so it converges. We put
$$x(t)=\lim_{n \rightarrow +\infty} x_n(t) .$$ 
Then  we have
$$-1<x_0(t)\leq x(t)\leq 1.$$
But, the equality $g_t(x_n(t))=2n\pi$ entails $\lim_{n\rightarrow +\infty}\, g_t(x_n(t))=+\infty$. Whence, by continuity of $g_t$, we deduce that $x(t)=1$.
$\hfill\square$

Now, we proceed to the study of $\nu_t$. Observe that for $ z\in \mathbb{D}$,
$$\left|z\exp\left(  \frac{t(z+1)}{2(z-1)}\right)\right|=|z| \exp\left(  \frac{t}{2}\Re \left(\frac{z+1}{z-1}\right)\right)=|z| \exp \left(-\frac{t(1-|z|^2)}{2|z-1|^2}\right)<1,$$
so $\theta$ is an analytic self-map of $\mathbb{D}$. From the previous identity one can easily check that $|\theta_t(\xi)|=1$ if $\xi\in \mathbb{T}$. 
Thus $\theta_t$ is an inner function and hence $\nu_t$ is singular with respect to the Haar measure on the unit circle.
Moreover, $\nu_t$ coincides with the Aleksandrov-Clark measure associated to $\theta_t$ at 1 (cf. \cite{ross}). 
Define the Poisson transform of $\nu_t$ by
\begin{equation*}
P_{\nu_t}(z):=\Re[k_t(z)]=\int_{\mathbb{T}}\frac{1-r^2}{|\xi-z|^2}d\nu_t(\xi),\quad t>0.
\end{equation*}
Let $z=r\zeta$, where $r\in(0,1)$ and $\zeta\in\mathbb{T}$, then the above expression becomes
\begin{equation*}
P_{\nu_t}(r\zeta)=\int_{\mathbb{T}}\frac{1-r^2}{|\xi-r\zeta|^2}d\nu_t(\xi)=\frac{1-|\theta_t(r\zeta)|^2}{|1-\theta_t(r\zeta)|^2},\quad t>0.
\end{equation*}
The rightmost term approaches zero as $r\rightarrow 1$ except for the set of solutions to the equation $\theta_t(\zeta)=1$. By Lemma \ref{discret}, this set is equal to
 $(\zeta_n)_{n\in\mathbb{Z}}$. Which says that $\nu_t$ is placing no mass on $\mathbb{T}\setminus\{(\zeta_n)_{n\in\mathbb{Z}}\}$ and so (cf. \cite[Proposition 4.14]{ross})
$$\nu_t=\sum_{n\in\mathbb{Z}}c_n(t)\delta_{\zeta_n},\quad t>0$$
where $c_n(t)$ are positive constants given by
$$c_n(t)=\frac{1}{|\theta_t'(\zeta_n)|}=\frac{|\zeta_n-1|^2}{|\overline{\zeta_n}(\zeta_n-1)^2-t|}.$$
Then by writing $\zeta_n=e^{i\alpha_n}$, we obtain
$$c_n(t)=\frac{2(1-\cos \alpha_n)}{t+2(1-\cos\alpha_n)}.$$
So we have
$$\nu_t=\sum_{n\in\mathbb{Z}} \frac{2(1-\cos \alpha_n)}{t+2(1-\cos\alpha_n)} \delta_{e^{i\alpha_n}},\quad t>0.$$

\begin{rem}
Since $(x_n(t))_{n\in\mathbb{N}}$ is increasing, then one can see that
\begin{equation*}
{\rm supp}\ \nu_t\subset\left\lbrace e^{i\theta}\left|\ -\arccos (x_0(t))\leq \theta \leq \arccos (x_0(t))\right. \right\rbrace 
\end{equation*}
 where $\arccos\in[0,\pi]$. 
\end{rem}

Next, we shall show that the above interval becomes the whole unit circle asymptotically as $t$ goes to infinity.
\begin{pro}
 One has
\begin{equation*}
 \lim_{t\rightarrow +\infty} x_0(t)=-1.
\end{equation*}
\end{pro}
 { \bf Proof}
Since $\partial_x g_t$ does not vanishes for any fixed $(t,x)\in \mathbb{R}^*_+\times[-1,1]$, the implicit function theorem implies that there exists a 
continuously differentiable function $\varphi$ such that $x_0=\varphi(t)$ and
\begin{equation*}
 \varphi'(s)=-\frac{\partial_t g_t(s,x)}{\partial_x g_t(s,x)}=\frac{x^2-1}{s+1-x} \leq 0
\end{equation*}
for any $s$ in an open neighborhood of $t$ in $\mathbb{R}^*_+$. It follows that $x_0$ is decreasing with respect to $t$ and thus the limit of $x_0$ as 
$t\rightarrow +\infty$
exists. Let $l=\lim_{t\rightarrow +\infty}x_0(t)$, then $l\in[-1,1)$.
But since $g_t(x_0(t))=0$, then we get
\begin{equation*}
 \sqrt{\frac{1+x_0(t)}{1-x_0(t)}}=\frac{2\arccos(x_0(t))}{t}.
\end{equation*}
The LHS of the previous identity goes to zero as $t\rightarrow+\infty$. Thus the RHS must do the same.
This gives us $\sqrt{(1+l)/(1-l)}=0$, so we have $l=-1$.
$\hfill\square$

\subsection{Boolean stochastic calculus and boolean unitary Brownian motion}
\subsubsection{\bf Boolean stochastic calculus}
We shall recall the relevant facts from boolean stochastic calculus as expounded in \cite{benghorbel+surmann04}.
For a given Hilbert space $H$ denote by $\Gamma(H)=\mathbb{C}\oplus H$ the boolean Fock space over $H$. The vector $\Omega=\binom{1}{0}$ is called the vacuum vector.
The vacuum expectation $\Phi_{\Omega}:B(\Gamma(H))\rightarrow\mathbb{C}$ is the state defined by
$$\Phi_{\Omega}(A):=\langle\Omega,A\Omega\rangle_{\Gamma(H)}.$$
Let $U\in H$, the creation and annihilation operators are bounded operators on $\Gamma(H)$ defined as
$$A^*(U)\left( \begin{array}{c}
\alpha\\ 
V \\ 
\end{array}\right)=\left( \begin{array}{c}
0\\ 
 \alpha\, U\\ 
\end{array}\right),\quad 
A(U)\left( \begin{array}{c}
\alpha\\ 
V \\ 
\end{array}\right)=\left( \begin{array}{c}
\langle U,V\rangle_H\\ 
 0\\ 
\end{array}\right)$$
where $\alpha\in\mathbb{C}$ and $V\in H$. We have $A(U)$ and $A^*(U)$ are adjoint to each other on $\Gamma(H)$ and
$$A(U)\,A^*(V)=\langle U,V\rangle_H\, P_{\Omega},$$
where $P_{\Omega}$ is the projection on $\mathbb{C}\Omega$.

Now let $h$ be a separable Hilbert space and let $H=L^2(\mathbb{R}_+;h)$ defined by
$$L^2(\mathbb{R}_+;h)=\left\lbrace f:\mathbb{R}_+\rightarrow h\ \text{measurable}\ ;\int_0^\infty||f(t)||^2_h \,dt<\infty \right\rbrace .$$
For a given unit vector $u\in h$, let $u^t:=\chi_{[0,t[}\otimes u\in L^2(\mathbb{R}_+;h)$ defined by $u^t(s)=\chi_{[0,t[}(s) u\in h$ and call
$A_t^*(u):=A^*(u^t),A_t(u):=A(u^t)$. A bounded operator $F$ on $\Gamma(L^2(\mathbb{R}_+;h))$ is given by
\begin{equation*}
F=\left( 
\begin{array}{cc}
\alpha & \langle\gamma| \\ 
|\delta\rangle & r \\ 
\end{array}
\right)
\end{equation*}
where \begin{itemize}
\item[\textbullet] $\alpha \in \mathbb{C}$,
\item[\textbullet] $\langle\gamma| \in B(L^2(\mathbb{R}_+;h);\mathbb{C})$ with $\gamma \in L^2(\mathbb{R}_+;h)$ given by $\langle\gamma| (f):= \langle \gamma,f\rangle_{L^2(\mathbb{R}_+;h)} $,
\item[\textbullet] $|\delta\rangle\in B(\mathbb{C};L^2(\mathbb{R}_+;h))$ with $\delta \in L^2(\mathbb{R}_+;h)$ given by $\langle\delta| (\lambda):= \lambda \delta$,
\item[\textbullet] $r\in B(L^2(\mathbb{R}_+;h))$.
\end{itemize}
\begin{defi}
\begin{enumerate}
 \item A family of bounded operators $F=(F_t)_{t\ge 0}$ on $  \Gamma(L^2(\mathbb{R}_+;h))$ is an adapted process if for any $t>0$ 
$$F_t\left( \begin{array}{c}
\lambda\\ 
f \\ 
\end{array}\right)
=F_t\left( \begin{array}{c}
\lambda\\ 
f \chi_{[0,t[}\\ 
\end{array}\right),\quad \text{for all}\ \left( \begin{array}{c}
\lambda\\ 
f \\ 
\end{array}\right)\in \Gamma(H).$$
We denote by $\mathcal{L}(h)$ the subspace of $\Gamma(H)$ of all adapted processes.
\item An element $F=(F_t)_{t\ge 0}\in \mathcal{L}(h)$ is said to be locally square integrable process if  the map $t\mapsto F_t$ from $\mathbb{R}_+$ to
$B(\Gamma(H))$ is measurable and satisfies
\begin{itemize}
 \item $t\mapsto \alpha_t\in L^2_{loc}(\mathbb{R}_+;\mathbb{C})$,
\item $t\mapsto \gamma_t\in L^2_{loc}(\mathbb{R}_+;L^2(\mathbb{R}_+;h))$,
\item $t\mapsto \delta_t\in L^2_{loc}(\mathbb{R}_+;L^2(\mathbb{R}_+;h))$.
\end{itemize}
We denote by $\mathcal{L}^2(h)$ the subspace of $\mathcal{L}(h)$ of locally square integrable processes.
\end{enumerate}
\end{defi}

According to \cite{benghorbel+surmann04}, for any locally square integrable process $F$, one can define stochastic integrals $\int_0^tF_s\,dA_s(u)$ and
$\int_0^tdA_s^*(u)\,F_s$ by using It\^o-Riemann sums. Furthermore, the processes
$t\mapsto \int_0^tF_s\,dA_s(u)\ \text{and}\ t\mapsto \int_0^tdA_s(u)^*\,F_s$
are adapted and continuous in norm. 

\subsubsection{\bf Fock space realization of boolean unitary Brownian motion}

We shall now provide a realization of the boolean unitary Brownian motion on the Boolean Fock space. According to \cite{benghorbel+surmann04} there exists a unique family of unitary operators $(U_t)_{t\in\mathbb{R}_+}$, solution to the stochastic
differential equation
\begin{equation}\label{QSDE2}
 dU_t=(dA_t(u)-dA_t^*(u)-\frac{1}{2}P_{\Omega}dt)U_t,\quad U_0=1
\end{equation}
given by
\begin{equation*}
U_t=\left( 
\begin{array}{cc}
\exp(\frac{-t}{2}) & L_{t,-u,\delta} \\ 
L^*_{t,u,\gamma_1} & M^*_{t,u,\zeta}+id_{L^2(\mathbb{R}_+;h)} \\ 
\end{array}
\right)
\end{equation*}


where \begin{itemize}
\item[\textbullet]
$\gamma_1: s\mapsto \exp(-s/2),$
\item[\textbullet]$\delta_t(.)= \exp(-\frac{1}{2}(t-.)),$
\item[\textbullet]$\zeta_t(.)= -\exp(-\frac{1}{2}(t-.))\chi_{[0,t]}(.)$
\end{itemize}
and
\begin{itemize}
\item[\textbullet]
 $L_{t,u,\gamma_1}:L^2(\mathbb{R}_+;h)\rightarrow\mathbb{C},\quad f\mapsto\int_0^t<u,f(s)>\gamma_1(s)ds,$
\item[\textbullet] $L^*_{t,u,\delta}:\mathbb{C}\rightarrow L^2(\mathbb{R}_+;h),\quad \lambda\mapsto\lambda\delta_t(.)\chi_{[0,t)}(.)u,$
\item[\textbullet] $M^*_{t,u,\zeta}:L^2(\mathbb{R}_+;h)\rightarrow L^2(\mathbb{R}_+;h),\quad 
M^*_{t,u,\zeta}(f)(s)= \langle \zeta_s,f\rangle_{L^2(\mathbb{R}_+;h)} \chi_{[0,t)}(s)\,u.$
\end{itemize}

For  $i\geq1$, we put the following notations
\begin{align*}
 \zeta^i_t(.)= (-1)^i\frac{(t-.)^{i-1}}{(i-1)!}\exp\left(-\frac{1}{2}(t-.)\right)\chi_{[0,t)}(.)
\end{align*}
and
\begin{align*}
\delta^i_t(.)= (-1)^{i-1}\frac{(t-.)^{i-1}}{(i-1)!}\exp\left(-\frac{1}{2}(t-.)\right)
\end{align*}
where we put $\zeta^1_t(.)=\zeta_t(.)$ and $\delta^1_t(.)=\delta_t(.)$. 

We need the following results.
\begin{lem}\label{LEM} For any $i\geq1$,
 \begin{align*}
   &(M^*_{t,u,\zeta})^{i}= M^*_{t,u,\zeta^i},\\
 &L_{t,-u,\delta}\circ M^*_{t,u,\zeta^i}=L_{t,-u,\delta^{i+1}},\\
 &L_{t,-u,\delta^{i+1}} \circ L^*_{t,u,\gamma_1}=\frac{(-t)^{i+1}}{(i+1)!} \exp (-t/2).
\end{align*}
\end{lem}
{ \bf Proof} Let $f\in L^2(\mathbb{R}_+;h)$.
We prove the first equality by induction on $i$. The case $i=1$ is obvious. Assume the result hold for $i$, then we have
\begin{align*}
 &(M^*_{t,u,\zeta})^{i+1}(f)
=M^*_{t,u,\zeta}\circ M^*_{t,u,\zeta^i}(f)
\\&=\int_0^. \exp\left(\frac{-1}{2}(.-s)\right)\int_0^s \langle u,f(r)\rangle(-1)^{i}\frac{(s-r)^{i-1}}{(i-1)!}\exp\left(-\frac{1}{2}(s-r)\right)dr \chi_{[0,t)}(s)ds\chi_{[0,t)}(.)\, u
\\&=\int_0^. \langle u,f(r)\rangle\exp\left(\frac{-1}{2}(.-r)\right)\int_r^. (-1)^{i}\frac{(s-r)^{i-1}}{(i-1)!} \chi_{[0,t)}(s) ds dr \chi_{[0,t)}(.)\, u
\\&=\int_0^. \langle u,f(r)\rangle\exp\left(\frac{-1}{2}(.-r)\right) (-1)^{i+1}\frac{(.-r)^{i}}{i!} dr \chi_{[0,t)}(.)\, u
\\&= M^*_{t,u,\zeta^{i+1}}(f).
\end{align*}

For the second assertion we have
\begin{align*}
 &L_{t,-u,\delta}\circ M^*_{t,u,\zeta^i}(f)
\\&=-\int_0^t \exp\left(\frac{-1}{2}(t-s)\right)\int_0^s \langle u,f(r)\rangle(-1)^{i}\frac{(s-r)^{i-1}}{(i-1)!}\exp\left(-\frac{1}{2}(s-r)\right)dr \chi_{[0,t)}(s)ds
\\&=-\int_0^t \langle u,f(r)\rangle \exp\left(\frac{-1}{2}(t-r)\right)\int_r^t (-1)^{i}\frac{(s-r)^{i-1}}{(i-1)!} \chi_{[0,t)}(s)dsdr
\\&=-\int_0^t \langle u,f(r)\rangle \exp\left(\frac{-1}{2}(t-r)\right) (-1)^{i}\frac{(t-r)^{i}}{i!} dr
\\&= L_{t,-u,\delta^{i+1}}(f)
\end{align*}
For the last assertion obseve that 
\begin{align*}
 L_{t,-u,\delta^{i+1}} \circ L^*_{t,u,\gamma_1}&=-\int_0^t\delta_t^{i+1}(s)\gamma_{1}(s)ds
\end{align*}
Then the equality follows from a simple computation.
$\hfill\square$

Now we come to the basic theorem.
\begin{theo}\label{theo2}
 The process $(U_t)_{t\in\mathbb{R}_+}$ is a boolean unitary Brownian motion.
\end{theo}
{ \bf Proof} We have to show that the process $(U_t)_{t\in\mathbb{R}_+}$ has boolean independent increments and that these increments  $ U_{st}:=U_tU_s^{\star}$ have a stationary distribution. We start with the question of boolean independent increments.
For a given $u\in h$, let $dX_t(u)=-i(dA_t(u)-dA_t^*(u))$ then from (\ref{QSDE2}) we obtain
\begin{equation*}
 U_t-U_s=i \int_s^tdX_r(u)U_r-\frac{1}{2}\int_s^tP_{\Omega}U_rdr,
\end{equation*}
so that,
\begin{equation}\label{increments}
 U_{st}=1+i \int_s^tdX_r(u)U_rU_s^{\star}-\frac{1}{2}\int_s^tP_{\Omega}U_rU_s^{\star}dr.
\end{equation}
Thus the operator $U_{st}-1$ lives on the time interval $[s, t)$. 
But since operators to disjoint time intervals are boolean independent (c.f. \cite{benghorbel+surmann04}, Remark 2.2), we get the boolean independence property of the process $(U_t)_{t\in\mathbb{R}_+}$.
It remains to show that $U_{st}$ has a stationary distribution. 
>From the equality (\ref{increments}) we have $t\mapsto U_{st}$ is obtained by solving the same equation as $U$ but with respect to the boolean Brownian motion $(X_t(u)-X_s(u))_{t\geq s}$. 
This implies that the increments of $U$ have a stationary distribution. So it is enough 
 to prove that the distribution of $U_t$ is $\nu_t$ for all $t\geq 0$.
For $n\geq1$, we put
\begin{equation*}
P_{\Omega}U_t^n=\left( 
\begin{array}{cc}
\Phi_{\Omega}(U_t^{n})& K_{t,n} \\ 
0 & 0 \\ 
\end{array}
\right).
\end{equation*}
Then
\begin{align*}
P_{\Omega}U_t^{n+1}&=P_{\Omega}U_t^{n}\left( 
\begin{array}{cc}
\exp (\frac{-t}{2}) & L_{t,-u,\delta} \\ 
L^*_{t,u,\gamma_1} & M^*_{t,u,\zeta}+id \\ 
\end{array}
\right)
\\&=\left( 
\begin{array}{cc}
\exp (\frac{-t}{2})\Phi_{\Omega}(U_t^{n})+ K_{t,n}\circ L^*_{t,u,\gamma_1} & \Phi_{\Omega}(U_t^{n}) L_{t,-u,\delta}+ K_{t,n}\circ (M^*_{t,u,\zeta}+id)\\ 
0 & 0 \\ 
\end{array}
\right).
\end{align*}
By a simple identification we obtain
\begin{equation}\label{Sys}
 \begin{cases}

  \Phi_{\Omega}(U_t^{n+1})=\exp (-t/2) \Phi_{\Omega}(U_t^{n})+ K_{t,n}\circ L^*_{t,u,\gamma_1}\\
K_{t,n+1}=K_{t,n}\circ (M^*_{t,u,\zeta}+id) + \Phi_{\Omega}(U_t^{n}) L_{t,-u,\delta}
 \end{cases}.
\end{equation}
The second equation of (\ref{Sys}) gives by induction
\begin{align*}
 K_{t,n+1}=\sum_{k=1}^{n+1}\Phi_{\Omega}(U_t^{n+1-k}) L_{t,-u,\delta}\circ (M^*_{t,u,\zeta}+id)^{k-1}.
\end{align*}
Which in turn leads to
\begin{align*}
 K_{t,n+1}=\sum_{k=1}^{n+1}\Phi_{\Omega}(U_t^{n+1-k}) \sum_{i=0}^{k-1}\binom{k-1}{i}L_{t,-u,\delta}\circ (M^*_{t,u,\zeta})^{i}.
\end{align*}
Then, by Lemma \ref{LEM} we obtain 
\begin{equation*}
 K_{t,n+1}=\sum_{k=1}^{n+1}\Phi_{\Omega}(U_t^{n+1-k}) \sum_{i=0}^{k-1}\binom{k-1}{i}L_{t,-u,\delta^{i+1}}.
\end{equation*}
Now, by combining the last assertion of Lemma \ref{LEM} and  (\ref{Sys}) we get 
\begin{align*}
 \Phi_{\Omega}(U_t^{n+1})&=\exp (-t/2)\Phi_{\Omega}(U_t^{n})+\sum_{k=1}^{n}\Phi_{\Omega}(U_t^{n-k}) \sum_{i=0}^{k-1}\binom{k-1}{i}L_{t,-u,\delta^{i+1}} \circ L^*_{t,u,\gamma_1}
\\&=\exp (-t/2)\left[ \Phi_{\Omega}(U_t^{n})+\sum_{k=0}^{n-1}\Phi_{\Omega}(U_t^{n-1-k}) \sum_{i=0}^{k}\binom{k}{i}\frac{(-t)^{i+1}}{(i+1)!}\right].
\end{align*}
We introduce now the generating function
\begin{equation*}
 \eta(t,z)=\sum_{n=1}^{\infty}\Phi_{\Omega}(U_t^{n})z^n,\quad |z|<1.
\end{equation*}
We obtain
\begin{align*}
 \eta(t,z)-z \exp (-t/2)&=\exp (-t/2)\left[ z\eta(t,z)+z^2( 1+\eta(t,z))\sum_{n=0}^{+\infty} \left( \sum_{i=0}^{n}\binom{n}{i}\frac{(-t)^{i+1}}{(i+1)!}\right) z^n\right]
\end{align*}
which simplifies to
\begin{align*}
 \eta(t,z)&=\exp (-t/2)z(1+\eta(t,z))\left[1 +\sum_{n=0}^{+\infty} \left( \sum_{i=0}^{n}\binom{n}{i}\frac{(-t)^{i+1}}{(i+1)!}\right) z^{n+1}\right].
\end{align*}
Now using the expansion
\begin{align*}
 \exp \left( \frac{tz}{z-1}\right) =1 +\sum_{n=1}^{+\infty} \left( \sum_{i=0}^{n-1}\binom{n-1}{i}\frac{(-t)^{i+1}}{(i+1)!}\right) z^n,
\end{align*}
we derive
\begin{align*}
 \eta(t,z)&=\exp (-t/2) z(1+\eta(t,z)) \exp\left( \frac{tz}{z-1}\right) 
\\&=z(1+\eta(t,z)) \exp \left(\frac{t(1+z)}{2(z-1)}\right) .
\end{align*}
Finally, we have for any $t>0,\ \eta(t,.)=\psi_{\nu_t}$ and the distribution of $U_t$ is $\nu_t$. This ends the proof of Theorem \ref{theo2}.
$\hfill\square$

\begin{center}
  \textbf{Acknowledgment }                      
 \end{center}

I would like to thank U. Franz, who has through his discussions greatly contributed to the present paper and  T. Hasebe for his remarks and advices. I also would like to thank an anonymous referee for helpful comments and suggestions.

\end{document}